\documentclass[12pt]{article}
\usepackage{amssymb}
\usepackage{amsmath}
\newtheorem{prop}{Proposition}[section]
\newtheorem{thm}{Theorem}
\newtheorem{Lemma}[prop]{Lemma}
\newtheorem{remark}[prop]{Remark}
\newtheorem{definition}[prop]{Definition}

\newtheorem{ex}[prop]{Example}

\newcommand{\fix}{\mbox{\rm Fix}}

\def\R{\Bbb R}   \def\Z{\Bbb Z} \def\T{\Bbb T}
\newcommand{\C}{{\Bbb C}}  \def\T{\Bbb T}

\begin{document}
\title{Free $\Z^p$-actions on the three dimensional torus.}
\author{Eduardo Fierro Morales\footnote{Partially supported by DGIP of the Universidad Cat\'olica
del Norte, Antofagasta, Chile} \, and Richard Urz\'{u}a-Luz \footnote{ Partially
supported by Fondecyt \# 1100832 and DGIP of the Universidad Cat\'olica del
Norte, Antofagasta, Chile}}

\date{}
\maketitle

\begin{abstract}

We show that for each natural $p\geq 2$, \ the Lefschetz fixed point theorem is optimal when applied to $\Z^{p}$-actions by homeomorphisms on the
three dimensional torus $\T^3$. More precisely, we show that for a spectrally unitary  $\Z^p$-action ${\bf A}$ on the first
homology group $H_1(\T^3,\Z )$ with trivial fixed point set, there exists a free $\Z^p$-action by real analytic
diffeomorphisms of $\T^3$ whose induced $\Z^p$-action on $H_1(\T^3,\Z)$ is the action ${\bf A}$. In particular, we establish the normal form for this type of actions.
\end{abstract}

\section{\textbf{Introduction.}}
Let $p,q$ be natural numbers. By a $\Z^p$-action of class $C^r$(resp. analytic) on the torus $\T^q=\R^q/\Z^q$ we mean a homomorphism $\varphi$ of $\Z^p$ into Diff$^r\left( \T^{q}\right)$, where Diff$^r\left(\T^q\right)$ is the group of $C^r$ diffeomorphisms of the torus $\T^q$, $0\leq r\leq \infty$ (resp. analytic). It is well known that the $\Z^p$-actions $\varphi$ of class $C^r$ on the torus $\T^q$  are closely related to the $\Z^p$-actions  on $H_1(\T^q,\Z)= \Z^q$. In fact, any  $\Z^p$-action of class $C^r$ on the torus $\T^q$ induces a $\Z^p$-action $\varphi _*$  on  $H_1(\T^q,\Z)$, which acts by automorphisms of $\Z^q$, called the {\it induced action } by $\varphi$. When $\varphi$ is free, that is $\varphi(\ell)$ has no fixed points, for each $\ell \in \Z^p\smallsetminus\{0\}$, it follows from the Lefschetz fixed point theorem that $1$ is an eigenvalue of $\varphi_*(\ell)$ for each $\ell\in\Z^p$.

We say that a  $\Z^{p}$-action ${\bf A}$ on $H_1(\T^q,\Z)$ is {\it spectrally unitary} if ${\bf A}$ acts by automorphisms of $\Z^q$ and $1$ belongs to the spectrum of ${\bf A}(\ell )$ for each $\ell \in \Z^p$. Notice that any free $\Z^p$-action $\varphi$ on the torus $\T^q$ induces a spectrally unitary $\Z^p$-action on $\Z^q$.

We will denote by ${\rm Fix}({\bf A})$ the set of fixed points of ${\bf A}$, i.e.,
\begin{equation}
{\rm Fix}({\bf A})=\{k\in \Z^q|{\bf A}(\ell )k=k \text{ for all }\ell \in \Z^p\}.
\end{equation}

In this paper, we are interested in the study of the existence of finite orbits for $\Z^p$-actions on the torus $\T^q$.
This study is related to the classification of the cohomologically  rigid $\Z^p$-actions on manifolds, in the sense that the cohomology group of the action is isomorphic to de Rham cohomology of the torus $\T^p$, see~\cite{P-S}. In fact, for $p=1$, the cohomologically  rigid diffeomorphisms of the torus $T^q$, where $q\leq3$, are $C^\infty$-conjugate to Diophantine translations ~\cite{L-S}. Recently Forni, Flaminio, and Rodriguez-Hertz ~\cite{F-F-RH} proved that the only cohomologically  rigid affine transformations of a homogenous space are conjugated to toral diophantine translations.

For $p\geq2$, dos Santos in ~\cite{S} proved that the cohomologically rigid $\Z^p$-actions on $\T^q$, $1\leq q\leq2$ are smooth conjugations of affine cohomologically rigid actions. In order to extend this result to actions on $\T^3$, one difficulty is knowing whether the induced homology of these actions has an invariant direction. More precisely,  in ~\cite[Problem 2.8.]{S}, dos Santos poses the following {\it stronger question}: let $\varphi:\Z^p\rightarrow $Diff$^r\left(\T^q\right)$, $r\geq 0$, be an action and $\varphi _*$ be its induced action on $H_1(\T^q,\Z)$ such that ${\rm Fix}({\bf \varphi _*})=\{0\}$. Does $\varphi$ have a finite orbit?

 When $q=1$  this question has an affirmative answer, this follows from Poincar\'e classification for homeomorphisms of the circle ~\cite{K-H}. The question also has an affirmative answer for $\Z^p$-actions of class $C^1$ on the torus $\T^2$, which follows from Franks, Handel, and Parwani's work ~\cite{F-H-P}. In ~\cite{S} the author observed that every affine $\Z^p$-action $\varphi$ on the torus $T^q$ satisfying $\fix(\varphi_*)=\{0\}$ has a finite orbit. We recall that a  $\Z^p$-action $\varphi$ on the torus $\T^q$ is {\it affine} if $\varphi(\ell)$ is an affine transformation of the torus $\T^q$ for each $\ell\in\Z^p$.

In this paper, we give a negative answer to the {\it stronger question} for $q=3$, which we establish in the following theorem.

\begin{thm}\label{principal}
Let $p\geq 2$ and let ${\bf A}$ be a spectrally unitary $\Z^p$-action on $H_1(\T^3,\Z)$. Then there is a free real analytic
$\Z^p$-action $\varphi$ on the torus $\T^3$ whose induced action on the first homology group of $\T^3$ is ${\bf A}$.
\end{thm}

Any $\Z^p$-action ${\bf A}$ on $\Z^q$ that acts by automorphisms of $\Z^q$, induce an action $A: \Z^p \rightarrow {\rm GL}(q,\Z)$. Conversely, any action $A : \Z^p \rightarrow {\rm GL}(q,\Z)$ defines an action by automorphisms. The image of ${\bf A}$ is by definition the
abelian group of automorphisms of $\Z^q$ generated by ${\bf A}(\ell)$, where $\ell\in\Z^p$ and is isomorphic to the subgroup of matrices $A(\Z^p)\subset{\rm GL}(q,\Z)$.

To prove Theorem \ref{principal}, we developed the normal form for the spectrally unitary $\Z^p$-actions on $H_1(\T^3,\Z)$ with trivial fixed point set, which is established in the following theorem.

\begin{thm} \label{normal}
Let $p\geq 2$  and let ${\bf A}$ be a spectrally unitary $\Z^p$-action on $\Z^3$ satisfying ${\rm Fix}({\bf A})=\{0\}$.
Then, there exists a basis $\{w_1,...,w_p\}$ of $\Z^p$ such that after a conjugation by an automorphisms of $\Z^3$, ${\bf A}$ has the following normal form in the canonical basis of $\Z^3$
$$ A(w_1)=\left[
\begin{array}{ccc}
1&a &b  \\
0&-1  &0  \\
0&0  &-1
\end{array}\right],\; A(w_2)=\left[
\begin{array}{ccc}
-1&0  &c \\
0& -1 &d \\
0&0  &1
\end{array}
\right]
$$
where $ad+2(b+c)=0$ and $A(\ell)$ is the identity matrix for all $\ell=\sum^p _{i=1}\ell _i w_i\in \Z^p$ satisfying $\ell _1=\ell
_2=0$. In particular, the image of ${\bf A}$ is  isomorphic to the Klein four-group.
\end{thm}

When ${\rm Fix}({\bf A})\neq\{ 0\}$ Theorem \ref{principal} follows from  ~\cite[Theorem 3]{L1},  whose proof is based on the classification of the unipotent $\Z^p$-actions on $\Z^3$. In each of these cases we constructed a free affine $\Z^p$-action.

For $p\geq3$ and $q\geq4$, not every unipotent action admits a free affine action. For  $p=2$, we do not know.

The rest of this paper is organized as follows. In section 2, we characterize, under certain hypotheses, the abelian subgroups of ${\rm GL}(q,\Z)$ with two generators (Lemma 2.2 and Lemma 2.5). Using these results, we obtain the normal form for spectrally unitary $\Z^p$-actions on $H_1(\T^q,\Z)$ (Theorem 2). In section 3, we establish the proof of our principal result (Theorem 1). In section 4, we reformulate the {\it stronger question} for $\Z^p$-actions on the torus $\T^q$ preserving the Haar measure.

\section{\textbf{Spectrally unitary $\Z^p$- actions on $\Z^3$ with trivial fixed point set.}}
In this section we will prove Theorem \ref{normal}. For this, we need to establish some results for spectrally unitary
$\Z^p$-actions on $\Z^3$ with trivial fixed point set.
\begin{definition}
Let ${\bf A}$ and ${\bf B}$ be two $\Z^p$-actions on $\Z^q$ that acts by automorphisms of $\Z^q$ . We say that ${\bf A}$ is  {\it conjugate} to ${\bf B}$ if there is an automorphism $h$ of $\Z^q$ such that ${\bf B}(\ell)\circ h=h\circ{\bf A}(\ell)$ for all $\ell \in \Z^p$.
\end{definition}

Let $A$ and $B$: $\Z^p\rightarrow {\rm GL}(q,\Z)$ be two actions. We say that $A$ and $B$ are conjugate if there is $P\in {\rm GL}(q,\Z)$ such that $B(\ell)P=PA(\ell)$ for all  $\ell \in \Z^p$. Two actions ${\bf A}$ and ${\bf B}$ are conjugate via automorphism if and only if the corresponding actions $A$ and $B$ are conjugate via an element of ${\rm GL}(q,\Z)$.

A fundamental example in the proof of our main result is the $\Z^2$-action ${\bf A}$ on $\Z^3$ generated, in the canonical basis $\{e_1,e_2\}$ of $\Z^2$, by the matrices
$$
A(e_1)=\left[
\begin{array}{ccc}
1&a &b  \\
0&-1  &0  \\
0&0  &-1
\end{array}
\right]
 \mbox{ and }
 A(e_2)=\left[
\begin{array}{ccc}
-1&0  &c \\
0& -1 &d \\
0&0  &1
\end{array}
\right]
$$
where $ad+2(b+c)=0$.
Since ${\bf A}(2e_1)={\bf A}(2e_2)=Id$, it follows that ${\bf A}$ is a spectrally unitary and an easy computation shows that ${\rm Fix}({\bf A})=\{0\}$. In fact, the Theorem 2 states that, for an adequate basis of $\Z^p$; any $\Z^p$-action on $\Z^3$ is essentially the action ${\bf A}$. For this, we need to show the following lemmas.

\begin{Lemma}\label{forma normal}
Let $N$ and $M$ be in ${\rm GL}(3,\Z)$ satisfying the following properties:
\begin{enumerate}
\item [i.]$NM=MN$
\item [ii.]$1$ is an eigenvalue of $N^n\dot M^m$ for all $(n,m)\in \Z^2$.
\item [iii.] $\{k\in \Z^3/ N^n\dot M^m k=k, \mbox{ for all } (n,m)\in \Z^2 \}=\{0\}$.
\end{enumerate}
Then after conjugation by an element of ${\rm GL}(3,\Z)$, $N$ and $M$  have the following form
\begin{equation}
N=\left[
\begin{array}{ccc}
1&a &b  \\
0&-1  &0  \\
0&0  &-1
\end{array}
\right]
 \mbox{, }
 M=\left[
\begin{array}{ccc}
-1&0  &c \\
0& -1 &d \\
0&0  &1
\end{array}
\right] \label{nor}
\end{equation}
 In particular, the subgroup of ${\rm GL}(3,\Z)$ generated by $N$ and $M$ is isomorphic to the Klein four-group.
\end{Lemma}

\begin{remark}
 Lemma \ref{forma normal} shows that the image of a spectrally unitary $\Z^2$-action on $\Z^3$ with trivial fixed point set is isomorphic to the Klein four-group. The following example shows that this is not true for a general spectrally unitary $\Z^2$-action on $\Z^q$ with trivial fixed point set.
\end{remark}

\begin{ex}\label{klein}
Let $\bf{A}$ be the action generated in the canonical basis $\{e_1,e_2\}$ of $\Z^2$ by
$$
A(e_1)=\left[
\begin{array}{cccc}
1 & a & 0 & 0 \\
0 & -1 & 0 & 0 \\
0 & 0 & 0 & -1 \\
0 & 0 & 1 & -2
\end{array}
\right] , A(e_2)=\left[
\begin{array}{cccc}
-1 & 0 & 0 & 0 \\
0 & -1 & 0 & 0 \\
0 & 0 & 1 & 0 \\
0 & 0 & 0 & 1
\end{array}
\right]
$$
with $a\in \Z$. It is easy to see that  $\bf{A}$ is a spectrally unitary $\Z^2$-action on $\Z^4$, satisfying the hypothesis of Lemma \ref{forma normal}. Notice that  the cyclic subgroup generated by ${\bf A}(e_1)$ is not finite. In particular, the image of ${\bf A}$ is not isomorphic to the Klein four-group.
\end{ex}
{\bf Proof of Lemma \ref{forma normal}}. Let $N$ and $M$ be matrices in ${\rm GL}(3,\C)$ satisfying the hypothesis i, ii, and iii. The subspaces ${\rm Fix}(N)$, ${\rm Fix}(M)$, and ${\rm Fix}(MN)$ are non trivial
subspaces (because of condition ii) whose pairwise intersections are reduced to $\{0\}$ (because of conditions i and iii).
The direct sum ${\rm Fix}(N)\oplus {\rm Fix}(M)$ has dimension of at least two and intersects trivially ${\rm Fix}(NM)$ because
each subspace ${\rm Fix}(N)$ and ${\rm Fix}(M)$ is invariant by $NM$ (because of i). We can deduce that ${\rm Fix}(MN)$ has dimension one, that  ${\rm Fix}(N)\oplus {\rm Fix}(M)$ has dimension two (which of course implies that each space ${\rm Fix}(N)$, ${\rm Fix}(M)$ has dimension one), and

$$
\C^3={\rm Fix}(N)\oplus {\rm Fix} (M)\oplus {\rm Fix}(NM).
$$
Each of these three subspaces being invariant by both $M$ and $N$, these two matrices are simultaneously diagonalizable.

Let $a(N)$, $b(N)$, $c(N)$, $a(M)$, $b(M)$, and $c(M)$ be the diagonal coefficients of $N$ and $M$, respectively.  We can suppose that $a(N)=1$, $b(N)\not=1$, $c(N)\not=1$, $a(M)\not=1$, $b(M)=1$, and $c(M)=c(N)^{-1}\not=1$. We observe that for every pair of integer $n$ and $m$, one of the numbers $b(N)^n$, $a(M)^m$ and $c(N)^nc(M)^m=c(N)^{n-m}$ is equal to 1. Setting $n=2$ and $m=1$, one gets $b(N)=-1$; setting $n=1$ and $m=2$, one gets $a(M)=-1$, setting $n=1$ and $m=-1$, one gets $c(N)=-1$. This implies that $M$ and $N$ generate the Klein four-group.

Now, we show that there exists a basis of $\Z^3$ such that in this basis,  $N$ and $M$, have the desired form. In fact, we can assume by
~\cite{N} that there exists a basis $\{e_1,e_2,e_3\}$ of $\Z^3$ such that

$$
N=\left[
\begin{array}{ccc}
1& \tilde{a}_{12}&\tilde{a}_{13}  \\
0& -1 &\tilde{a}_{23}  \\
0&0 & -1
\end{array}
\right]
$$
In particular, $\tilde{a}_{23}=0$ since $N^2=I$, where $I$ is the identity matrix. Considering that the eigenspace of $N$ associated with the eigenvalue $1$ has
dimension one and  $NMe_1=Me_1$, we have $Me_1=me_1$, and therefore $m$ is an eigenvalue of $M$. Hence, $m=\pm 1$. Thus from hypothesis
$(iii)$ we have that $m=-1$. Then, on this basis, $M$ has the form
$$
M=\left[
\begin{array}{ccc}
-1& \tilde{b}_{12}&\tilde{a}_{13}  \\
0& \tilde{b}_{22} &\tilde{b}_{23}  \\
0&\tilde{b}_{32} &\tilde{b}_{33}
\end{array}
\right].
$$
Let $B=\left[
\begin{array}{cc}
\tilde{b}_{22} &\tilde{b}_{23}  \\
\tilde{b}_{32} &\tilde{b}_{33}
\end{array}
\right]$. Then $\det M=-\det B=1$. Hence $\det B=-1$. Since $-1$ is an eigenvalue of multiplicity two, we have that $-1$ is an eigenvalue of $B$. Then it is easy to see that there exists a matrix $U\in {\rm GL}(2,\Z)$ such that

\begin{equation}
U^{-1}BU=\left[
\begin{array}{cc}
-1 & b_{23}  \\
0 &b_{33}
\end{array}
\right]
\end{equation}
Therefore
\begin{equation}
\left[
\begin{array}{ccc}
1&0&0\\
0& U&
\end{array}
\right]^{-1}N\left[
\begin{array}{ccc}
1&0&0\\
0& U&
\end{array}
\right]=\left[
\begin{array}{ccc}
1& a_{12}&a_{13}  \\
0& -1 &0  \\
0&0 & -1
\end{array}
\right]
\end{equation}
and
\begin{equation}
\left[
\begin{array}{ccc}
1&0&0\\
0& U&
\end{array}
\right]^{-1}M\left[
\begin{array}{ccc}
1&0&0\\
0& U&
\end{array}
\right]=\left[
\begin{array}{ccc}
-1& b_{12}&b_{13}  \\
0& -1&b_{23}  \\
0&0&b_{33}
\end{array}
\right]
\end{equation}

Since $\det M=1$, we have  $b_{33}=1$. Now, from hypothesis $(i)$, we have $b_{12}=0$ and $a_{12}b_{23}=-2(a_{13}+b_{13})$, thus finishing the
proof.

 $\hfill\square$

An immediate consequence of  Lemma \ref{forma normal} is the following result.

\begin{Lemma}\label{maximal}
Let $N$ and $M$ be in ${\rm GL}(3,\Z)$ under the same hypothesis of Lemma \ref{forma normal} and $C \in {\rm GL}(3,\Z)$.
If $C$ commutes with $N$ and $M$, and  $1$ is an eigenvalue of $N^nM^mC^s$ for all $(n,m,s)\in \Z^3$, then $C$ belongs to the group
generated by $N$ and $M$.
\end{Lemma}
{\bf Proof.} By Lemma \ref{forma normal} we observe that the decomposition
$$
\C^3={\rm Fix} (N)\oplus {\rm Fix}(M)\oplus {\rm Fix} (NM)
$$
is invariant by every element $C$ in ${\rm GL}(3,\C)$ that commutes with $M$ and $N$. Writing $a(C)$, $b(C)$, and $c(C)$ the diagonal coefficients of $C$, we  can deduces that if the hypothesis of Lemma \ref{maximal} is satisfied, then for every integer $n$,$m$, $s$, one of the coefficient $(-1)^ma(C)^s$, $(-1)^nb(C)^s$, $(-1)^{n+m}c(C)^s$ is equal to 1. It is straight forward that $C$ is equal to $I$,  $M$, $N$, or $MN$.

$\hfill\square$

\begin{remark} Lemma \ref{maximal}  states that the abelian  group $\left<N,M\right>$ generated by $N$ and $M$ is "maximal" in ${\rm GL}(3,\Z)$ in the following sense: if ${\rm G}$ is an abelian subgroup of ${\rm GL}(3,\Z)$ that contains $N$ and $M$, and $1$ is an eigenvalue of each element of ${\rm G}$, then ${\rm G}=\left<N,M\right>$.
\end{remark}
Let $p$ be a natural number. For each $i=1,2,...,p$, we denote by ${\rm G}_i$ the subgroup of $\Z^p$ defined by ${\rm G}_i=\{(a_1,a_2,...,a_p)\in
\Z^p\mid a_i=0\}$.

Given a $\Z^p$-action ${\bf A}$ on $H_1(\T^3, \Z)$ acting by automorphisms of $\Z^3$,  we shall denote by ${\bf A}_i$ the restriction of the action ${\bf A}$ to the subgroup $G_i$. If $\{e_1,...,e_p\}$ is the canonical basis of $\Z^p$,  it is easy to see that ${\rm Fix}({\bf A}_i)=\cap _{j\neq i}{\rm Fix}({\bf A}(e_j))$, for all $i$. Thus we obtain the following lemma.
\begin{Lemma}\label{fijo}
Let ${\bf A}$ be an spectrally unitary $\Z^p$-action on $\Z^3$. If $p\geq 3$ and  ${\rm Fix}({\bf A})=\{0\}$, then ${\rm Fix}({\bf A}_i)=\{0\}$, for some $i$ such that $1\leq i\leq p$.

\end{Lemma}
{\bf Proof.} Suppose the lemma is false. Then there exists a set
of nonzero vectors $\{v_1,v_2,...,v_p\}\subset\Z^3$ such that $v_i\in{\rm Fix}({\bf A}_i)$ for all $i=1,2,...,p$.

\noindent Assume $\{v_1,v_2,...,v_p\}$ is linearly dependent. Then there exist a natural number $k$, $k\leq p$, and rational numbers
$\lambda  _1$, $\lambda _2$,...,$\lambda _p$ such that
$$
v_k=\sum _{i\neq k}\lambda _iv_i.
$$
It follows that ${\bf A}(e_k)v_k=v_k$ and $v_k\in{\rm Fix}({\bf A})$. Therefore ${\rm Fix}({\bf A})\neq\{0\}$, a contradiction with the hypothesis.

\noindent Now suppose $\{v_1,v_2,...,v_p\}$ is linearly independent. Then $p=3$ and there exist rational numbers $a_{11},a_{21}$, and $a_{31}$ such that ${\bf A}(e_1)v_1=a_{11}v_1+a_{21}v_2+a_{31}v_3$. Since ${\rm Fix}({\bf A})=\{0\}$
and ${\bf A}(e_1)v_1$ is invariant by ${\bf A}(e_k)$ for $k=2,3$, we have $a_{11}=-1$ and $a_{21}=a_{31}=0$. Therefore ${\bf A}(e_1)v_1=-v_1$. In a similar way we prove that ${\bf A}(e_2)v_2=-v_2$ and ${\bf A}(e_3)v_3=-v_3$. Thus, $1$ is not an eigenvalue of ${\bf  A}(e_1){\bf A}(e_2){\bf A}(e_3)$, which contradicts the hypothesis. $\hfill \square$
\medskip

\noindent {\bf Proof of Theorem \ref{normal}.} We proceed by induction on $p$. If $p=2$, we may assume $A(e_1)$ and $A(e_2)$ satisfy the hypothesis of Lemma \ref{forma normal}. Hence, after conjugation by an element of ${\rm GL}(3,\Z)$, we may assume that $ A(e_1)$ and $A(e_2)$ are in the normal form $(\ref{nor})$. Thus, the  canonical basis $\{e_1,e_2\}$ of $\Z^2$ has the required properties.

Assume the theorem holds for $p-1$; we will prove that it also holds for $p$. Since ${\rm Fix}({\bf A})=\{0\}$, by Lemma \ref{fijo}
we see that ${\rm Fix}({\bf A}_i)=\{0\}$ for some $i=1,2,...,p$. Without loss of generality, we can assume that $i=p$. Then by the  induction
hypothesis there exists a basis $\{v_1,...,v_{p-1}\}$ of $G_p$ such that  $ A(v_1)$ and $A(v_2)$ are in the normal form by $(\ref{nor})$
and $ A(v_j)=I$ for all $j>2$. Since  $ A(e_p)$  commutes with $A(v_1)$ and $A(v_2)$ then from Lemma \ref{maximal}
we have $A(e_p)= A(v)$ for some $v\in\{0,v_1,v_2,v_1+v_2\}$. So the basis $\{w_1,w_2,...,w_{p-1}, w_{p}\}$ of $\Z^p$,
where $w_k=v_k$ for all $ k < p$ and $w_p=v+e_p$, has the required properties, proving Theorem \ref{normal}.
$\hfill \square$

\section{\textbf{Proof of Theorem \ref{principal}}}
To prove Theorem \ref{principal}, we use the following proposition.
\begin{prop}\label{subgroup}
Let $G$ be a torsion-free group and $\varphi$ a $G$-action on a set $X$. Let $H$ be a normal subgroup of $G$ of finite index. Then $\varphi$ is a free $G$-action if and only if the restriction of the action $\varphi$ to $H$ is free.
\end{prop}
{\bf Proof.} If $\varphi $ is a free $G$-action then any restriction of $\varphi $ is free. Conversely, suppose there exists $g\in G$ such that $\varphi (g)$ has a fixed point $x\in X.$ Let $m=\left[ G:H\right]$. Then $g^m\in H$ and $\varphi (g^m)x=\varphi ^m(g)x=x$. Since the restriction of $\varphi $ to $H$ is free then $g^m=e,$ where $e$ is the identity element of $G$. Since $G$ is a torsion-free group, we have that $g=e$.
$\hfill \square$

Suppose ${\rm Fix}({\bf A})=\{0\}$.  By the normal form of the action ${\bf A}$ given in Theorem \ref{normal}, there exists a basis $\{w_1,...,w_p\}$ of $\Z^p$ such that $A(w_1)=N$, $ A(w_2)=M$, and $A(w_j)=I$ for $j>2$, where the matrices $N$ and $M$ are defined as in Lemma \ref{forma normal}. Consider the diffeomorphisms $\varphi _i$, $i=1,...,p$ of the torus $\T^3$ whose linear part is determined by the matrix $A(w_i)$, $i=1,...,p$ respectively, and which are given  on the covering $\R^3$ by
\begin{eqnarray}
\widetilde{\varphi} _1(x,y,z)&=&(x+ay+bz+f_1(z)+r, -y+g_1(z),-z)
\label{varhi1}\\
\widetilde{\varphi} _2(x,y,z)&=&(-x+cz+f_2(z),-y+dz+g_2(z),z+\tfrac{1}{2})
\label{varhi2}\\
\widetilde{\varphi} _j(x,y,z)&=&(x+f_j(z),y+g_j(z),z), \mbox{ for }j>2
\end{eqnarray}
where  $f_i$ and $g_j$, with $i,j=1,2,...,p$, are $\Z$- periodic functions of the variable $z$ with Lebesgue measure zero and $r$ is a
rational number. The family of diffeomorphisms $\widetilde{\varphi}_i$ defines  a $\Z^p$-action on the torus $\T^3$ if and only if
\begin{equation}
\widetilde{\varphi} _i\circ\widetilde{\varphi} _j-\widetilde{\varphi} _j\circ\widetilde{\varphi} _i \in \Z^3
\label{eq7}
\end{equation}
for all $i,j=1,2,...,p$.

Since the equations in $(\ref{eq7})$ are trivially satisfied for $i\geq 3$ and $j\geq 3$,  we need to verify the equations in $(\ref{eq7})$ for $i=1,2$ and $j=1,2,3,...,p$. Let $r=-\tfrac{b}{4}$,  $R(z)=-z$, and $T_{\tfrac{1}{2}}(z)=z+\tfrac{1}{2}$. First,  we consider the case $i=1$ and $j=2$, then we have
\begin{eqnarray*}
\widetilde{\varphi} _2\circ \widetilde{\varphi} _1( x,y,z)&=&\widetilde{\varphi} _2(\widetilde{\varphi} _1(x,y,z))
 \\&=&\widetilde{\varphi} _2( x+ay+bz+f_1(z)+r,-y+g_1(z),-z)\\
&=&(-x-ay-bz-f_1(z) -r-cz+f_2(-z),\\& &
y-g_1(z)-dz+g_2(-z),-z+\tfrac{1}{2})
\end{eqnarray*}

\begin{eqnarray*}
\widetilde{\varphi} _1\circ \widetilde{\varphi} _2( x,y,z)&=& \widetilde{\varphi} _1(\widetilde{\varphi} _2(x,y,z))
 \\&=&\widetilde{\varphi}_1(-x+cz+f_2(z),-y+dz+g_2(z),z+\tfrac{1}{2})\\&=&(-x+cz+f_2(z)-ay+adz+ag_2(z) \\& &
+bz+\tfrac{b}{2}+f_1(z+\tfrac{1}{2})+r,y
-dz-g_2(z)+\\& &g_1(z+\tfrac{1}{2}),-z-\tfrac{1}{2}).
\end{eqnarray*}

It is easy to see that $\widetilde{\varphi} _2$ and $\widetilde{\varphi} _1$ satisfy the equation $(\ref{eq7})$, if and only if
\begin{equation}
f_2\circ R-f_2=ag_2+f_1\circ T_{\tfrac{1}{2}}+f_1 \mbox{ and }
g_2+g_2\circ R=g_1+g_1\circ T_{\tfrac{1}{2}} \label{eqf12}.
\end{equation}
When  $i=1$ and $j\geq 3$ we have
\begin{eqnarray*}
\widetilde{\varphi} _j\circ \widetilde{\varphi} _1( x,y,z) &=& \widetilde{\varphi} _j(\widetilde{\varphi} _1(x,y,z))
\\&=&\widetilde{\varphi} _j( x+ay+bz+f_1(z)+r,-y+g_1(z),-z)\\&=& (x+ay+bz+f_1(z)+r+f_j(-z),\\& & -y+g_1(z)+g_j(-z),-z)
\end{eqnarray*}
and
\begin{eqnarray*}
\widetilde{\varphi} _1\circ \widetilde{\varphi} _j( x,y,z) &=&\widetilde{\varphi} _1(\widetilde{\varphi} _j(x,y,z))
\\&=&\widetilde{\varphi}_1(x+f_j(z),y+g_j(z),z)\\&=&(x+f_j(z)+ay+ag_j(z) +bz+f_1(z)+r,\\&&-y-g_j(z)+g_1(z),-z).
\end{eqnarray*}
Thus the condition  $(\ref{eq7})$ is verified  if
\begin{equation}
f_j\circ R-f_j=ag_j \mbox{ and } g_j\circ R=-g_j \label{eqf1j}.
\end{equation}
In the case where $i=2$ and $j\geq 3$ we have
\begin{eqnarray*}
\widetilde{\varphi} _j\circ \widetilde{\varphi} _2( x,y,z) &=& \widetilde{\varphi} _j(\widetilde{\varphi} _2(x,y,z)) \\&=&\widetilde{\varphi}_j(-x+cz+f_2(z),-y+dz+g_2(z),z+\tfrac{1}{2})\\&=&(x+cz+f_2(z)+f_j(z+\tfrac{1}{2}),-y+cz+g_2(z)\\&+&g_j(z+\tfrac{1}{2}),z+\tfrac{1}{2})
\end{eqnarray*}
and
\begin{eqnarray*}
\widetilde{\varphi} _2\circ \widetilde{\varphi} _j( x,y,z) &=&\widetilde{\varphi} _2(\widetilde{\varphi} _j(x,y,z))
 \\&=&\widetilde{\varphi}_2(x+f_j(z),y+g_j(z),z)\\&=&(-x-f_j(z)+cz+f_2(z),\\&
&-y-g_j(z)+dz+g_2(z),z+\tfrac{1}{2}).
\end{eqnarray*}
In this case, condition $(\ref{eq7})$ is satisfied if
 \begin{equation}
f_j\circ T_{\tfrac{1}{2}}=-f_j \mbox{ and } g_j\circ T_{\tfrac{1}{2}}=-g_j. \label{eqg2j}
\end{equation}
We choose $f_{1}\left( z\right) =\tfrac{\alpha _{1}}{2}\cos 2\pi z$, $g_{1}\left(z\right) =-\tfrac{\alpha _{1}}{2}\sin 2\pi z $, $f_{2}\left( z\right) =-\tfrac{\alpha _{2}}{2}\cos 2\pi z+\tfrac{a}{4}\alpha_{2}\sin 2\pi z$, $g_{2}\left( z\right) =-\tfrac{\alpha _{2}}{2}\sin 2\pi z$, and for $j>2$ we choose  $f_{j}\left( z\right) =\alpha _{j}(\cos 2\pi z-\tfrac{a}{2}\sin 2\pi z)$ and $g_{j}\left( z\right) =-\alpha _{j}\sin 2\pi z$, where $\alpha _1,\alpha _2,...,\alpha _p$ are arbitrary real numbers. A straightforward calculation shows that equations $(\ref{eqf12})$, $(\ref{eqf1j})$, and $(\ref{eqg2j})$ are verified.
  Thus the family of diffeomorphisms  $\widetilde{\varphi} _j$, with $1\leq j\leq p$, defines a  $\Z^p$-action on the torus $\T^3$.
In what follows, we assume that the real numbers $\alpha _1,\alpha _2,...,\alpha _p$ are  linearly independent over the field of algebraic numbers.
We will show that the $\Z^p$-action $\varphi$ defined  by the family $\widetilde{\varphi} _j$ is a  free $\Z^p$-action.
For this, by Proposition \ref{subgroup} it suffices to show that the restriction of $\varphi$ to the subgroup ${\rm H}=\{2\ell_1 w_1+2\ell_2 w_2+...+ \ell_p w_p| \ell_1,...,\ell_p\in \Z \}$ is free.
From $(\ref{varhi1})$ and $(\ref{varhi2})$ we have that
$$
\widetilde{\varphi} _1^2(x,y,z)=(x+\alpha _1\cos 2\pi z-\tfrac{a}{2}\alpha _1 \sin 2\pi z+2r,y+\alpha _1\sin 2\pi z, z)
$$
and
$$
\widetilde{\varphi} _2^2(x,y,z)=(x+\alpha _2\cos 2\pi z-\tfrac{a}{2}\alpha _2 \sin 2\pi z+\tfrac{c}{2},y+\alpha _2\sin 2\pi z+\tfrac{d}{2}, z+1),
$$
thus we see that
$$
\widetilde{\varphi} _1^{2\ell _1}\circ \widetilde{\varphi} _2^{2\ell _2} \circ...\circ\widetilde{\varphi}_p^{\ell _p}(x,y,z)=
$$
$$
(x+\left(\sum\nolimits^p_{j=1}\ell _j\alpha _j\right)\cos 2\pi z-\tfrac{a}{2}\left(\sum\nolimits^p_{j=1}\ell _j\alpha _j\right)\sin 2\pi z+2r+\ell _1\tfrac{c}{2},
$$
\begin{equation}
y+\left(\sum\nolimits^p_{j=1}\ell _j\alpha _j\right)\sin 2\pi z+\ell _2\tfrac{d}{2},z+\ell _2).
\label{eqpar}
\end{equation}
If
$$
\varphi(2\ell_1 w_1+2\ell_2 w_2+...+ \ell_p w_p)=\varphi_1^{2\ell _1}\circ\varphi_2^{2\ell _2} \circ...\circ\varphi_p^{\ell _p}
$$
has a fixed point $\pi(x,y,z)\in \T^3$, then, on the covering $\R^3$, there exists $(n_1,n_2,n_3)\in \Z^3$ such that
$$
\widetilde{\varphi} _1^{2\ell _1}\circ \widetilde{\varphi} _2^{2\ell _2}\circ...\circ\widetilde{\varphi}_p^{\ell _p}(x,y,z)=(x,y,z)+(n_1,n_2,n_3).
$$
From $(\ref{eqpar})$ we see that
\begin{equation}
\left(\sum\nolimits^p_{j=1}\ell_j\alpha _j\right)\sin 2\pi z=n _2-\ell _2\tfrac{d}{2}
 \label{eqsin}
\end{equation}
and
\begin{equation}
\left(\sum\nolimits^p_{j=1}\ell _j\alpha _j\right)\cos 2\pi z=-2\ell _1r-\ell _2\tfrac{c}{2}+\tfrac{a}{2}\left(\sum\nolimits^p_{j=1}\ell_j\alpha _j\right)\sin 2 \pi z+n _1. \label{eqcos}
\end{equation}
Replacing $(\ref{eqsin})$ in equation $(\ref{eqcos})$ we get
\begin{equation}
\left(\sum^p _{j=1}\ell_j\alpha _j\right)\cos 2\pi z=-2\ell _1r-\ell _2\tfrac{c}{2}+\tfrac{a}{2}n_2-a\ell _2\tfrac{d}{4}+n _1
\label{eqcos1}
\end{equation}

From $(\ref{eqsin})$ and $(\ref{eqcos1})$ we see that

$$
\left(\sum\nolimits^p_{j=1}\ell _j\alpha _j\right)^2=\left( n_2-\ell _2\tfrac{d}{2}\right)^2+\left(-2\ell _1r-\ell_2\tfrac{c}{2}+\tfrac{a}{2}n_2-a\ell _2\tfrac{d}{4}+n _1 \right)^2.
$$
Since the family $\{\alpha _j| j=1,2,...,p\}$ was chosen to be linearly independent over the field of algebraic numbers, we have $\ell _1=\ell_2=...=\ell _p=0$. Thus, Theorem \ref{principal} is proved.

$\hfill\square$

\section{Comments.}
Since the {\it strong question} is false in general, we can reformulate this question for $\Z^p$-actions on the torus $T^q$ that preserve the Haar measure and with its corresponding affine action conjugated  to its linear part.

Let $\varphi$ be a $\Z^p$-action on the torus $\T^q$  and ${\bf A}$ be the induced action on $H_1(\T^q,\Z)$. Consider for each $\ell\in\Z^p$ a lifting $\widetilde{\varphi}(\ell):\R^q\rightarrow \R^q$ on the covering $\pi:\R^q\rightarrow T^q$ of $\varphi(\ell):\T^q\rightarrow \T^q$. Thus $\widetilde{\varphi}(\ell)$ can  be written as
\begin{equation}
\widetilde{\varphi}(\ell)={\bf A}(\ell)+\widetilde{F}(\ell)
\label{ques}
\end{equation}
where $\widetilde{F}(\ell)=F(\ell)\circ\pi$ and $F(\ell):\T^q\rightarrow \R^q$.

Denote by  $a$ the  $\Z^p$-action on $\T^q$ which acts by automorphisms of $\T^q$, defined by the action ${\bf A}$, i.e.
$$
a\left( \ell \right) \circ \pi =\pi \circ {\bf A}\left( \ell \right)
$$
for all $\ell \in \Z^p.$

\bigskip
An application $\alpha :\Z^p\rightarrow \T^q$ is 1-cocycle over the $\Z^p$-action $a$ if for each $\ell$, $\ell ^{\prime }$ $\in\Z^p$

$$
\alpha \left( \ell +\ell ^{\prime }\right) =\alpha \left( \ell \right) \cdot
a\left( \ell \right) \alpha \left( \ell ^{\prime }\right)
$$

A 1-cocycle is trivial if there exist a point $x_0\in \T^q$ such that
$$
\alpha \left( \ell \right) =a\left( \ell \right) \left( x_0\right)
^{-1}x_0
$$
for all  $\ell \in $ $\Z^p.$
A leading example of 1-cocycle over $a$ is given by invariant probability measures of the action.  Indeed, if $\mu$ is an invariant probability measure for $\varphi$, then the application  $\alpha_{\mu}$ defined by
$$
\alpha _{\mu }\left( \ell \right) =\pi \left( \mu \left( F\left( \ell\right) \right) \right)
$$
for all $\ell \in $ $\Z^p$, where $F$ is as in  $(\ref{ques})$, is a 1-cocycle over $a$.

\bigskip
Now, let $\varphi$ be a $\Z^p$-action on  $\T^q$  that preserves the Haar measure $m$ of $\T^q$. Suppose that ${\rm Fix}({\bf A})=\{0\}$ and that the 1-cocycle $\alpha_{m}$ is trivial. We propose the following question:  does the $\varphi$ action has a finite orbit?

The theorem 1 gives a negative answer to the {\it stronger question} in ~\cite{S}, but the 1-cocycle induced by the Haar measure is not trivial.

\subsection*{\textbf{ Acknowledgments.}} The authors want to thank Juan Rivera Letelier and Nathan M. dos Santos for helpful discussions.
\section*{}

\bigskip
Richard Urz\'ua Luz
\\ Universidad Cat\'olica del Norte, \\ Casilla 1280, Antofagasta, Chile.
\\ rurzua@ucn.cl
\bigskip
 \\ Eduardo Fierro Morales
\\ Universidad Cat\'olica del Norte, \\ Casilla 1280, Antofagasta, Chile.
\\ efierro@ucn.cl
\end{document}